\documentclass[letterpaper]{amsart}[12pt]

\usepackage{array, tabularx, multirow}

\usepackage[mathscr]{eucal}

\usepackage[utf8x]{inputenc}
\usepackage[english]{babel}
\usepackage{mathrsfs}
\usepackage{comment}  
\usepackage{enumerate}
\usepackage{amsmath}
\usepackage{pifont}  
\usepackage{amsthm}
\usepackage{appendix}
\usepackage{amssymb, wasysym}
\usepackage{color}
\usepackage[dvipsnames]{xcolor}
\usepackage{hyperref} 
\usepackage{epigraph}
\usepackage{lipsum}
\usepackage[displaymath,mathlines, pagewise]{lineno}
\hypersetup{pdfpagemode={UseOutlines},
bookmarksopen=true,
bookmarksopenlevel=0,
hypertexnames=false,
colorlinks=true,
citecolor=black,
linkcolor=black,
urlcolor=magenta,
pdfstartview={FitV},
unicode,
breaklinks=true,
}
\usepackage{graphicx}

\newcommand{\rtt}{\upharpoonright} 
\newcommand{\y}{\ \wedge\ }   
\newcommand{\mm}{\setminus}

\renewcommand{\subset}{\subseteq}

\newcommand{\id}{\mathrm{id}}
\newcommand*\rfrac[2]{{}^{#1}\!/_{#2}} 
\newcommand{\nom}{n\in\omega}
\newcommand{\ENV}{\mathrm{Env}}
\newcommand{\sat}{\mathrm{sat}}

\newtheorem{theorem}{Theorem}[section]
\newtheorem{lemma}[theorem]{Lemma}

\newtheorem{proposition}[theorem]{Proposition}
\newtheorem{corollary}[theorem]{Corollary}
\theoremstyle{definition}
\newtheorem{mydef}[theorem]{Definition}
\newtheorem{obs}[theorem]{Observation}

\newtheorem{example}[theorem]{Example}
\newtheorem{preg}[theorem]{Question}

\def\cI{{\mathcal{I}}}
\def\cJ{{\mathcal{J}}}
\def\cF{{\mathcal{F}}}
\def\cP{{\mathscr{P}}} 
\def\cR{{\mathcal{R}}}
\def\cC{{\mathcal{C}}}
\def\cX{{\mathcal{X}}}

\title{Familias independientes generalizadas} 

\begin{document}

\title{Generalized independence}

\author{Fernando Hern\'andez-Hern\'andez}
\address{Universidad Michoacana de San Nicolás de Hidalgo, México}
\email{fernando.hernandez@umich.mx}

\author{Carlos L\'opez-Callejas}
\address{Posgrado Conjunto en Ciencias Matemáticas UNAM-UMSNH, México}
\thanks{The first author’s research has been supported by CONACyT, Scholarship 733921 and PAPIIT IN104419}
\email{carloscallejas.math@gmail.com}

\begin{abstract}
    We explore different generalizations of the classical concept of independent families on $\omega$ following the study initiated by Fisher and Montoya \cite{fischer2019higher}. We show that under $\diamondsuit^*(\kappa)$ we can get strongly independent families on $\kappa$ of size $2^\kappa$ and present an equivalence of $\mathsf{GCH}$ in terms of strongly independent families. We merge the two natural ways of generalizing independent families through a filter or an ideal and we focus on the $\cC$-independent families, where $\cC$ is the club filter.  Also we show a relationship between the existence of $\cJ$-independent families and the saturation of the ideal $\cJ$.
\end{abstract}



\maketitle

\section*{Introduction}\pagenumbering{arabic}

Independent families are objects with strong combinatorial properties. Since their appearance in \cite{hausdorff2008zwei}, these families have been related to many other objects, such as almost disjoint families, ultrafilters and ideals. See for example \cite{Geschke}.

Independent families are naturally defined over the set of non-negative integers $\omega $; however, it is not clear what their natural generalization to larger cardinals should be. An \emph{independent family} on $\omega$ is a family $\cI\subseteq\cP(\omega) $ such that if $ S, T \subseteq\cI $ are finite and disjoint subfamilies then $ \bigcap S\setminus \bigcup T $ is infinite (we call this set a finite Boolean combination from $\cI$). In other words, on $\omega $, a family is independent if all its finite Boolean combinations are infinite. When we move to the case of an arbitrary cardinal $\kappa $ the notion of independence could be generalized in at least two different ways: the first would be by allowing larger Boolean combinations, that is, not only finite Boolean combinations but also the ones of length less than or equal to $\lambda$ for some given $\lambda$ and the second way would be to ask that finite Boolean combinations not only have infinite cardinality (or cardinality $\kappa$) but that they fulfill some notion of \textit{largeness}. 

The first of these generalizations is what is normally known in the literature as strongly independent families and these have recently been studied by Vera Fischer and Diana Montoya in \cite{fischer2019higher}. In the first section we define these families, we justify the reason for considering Boolean combinations of length less than $\kappa$ and we give a characterization of the Continuum Hypothesis in terms of the existence of one of these families on $\omega_1$, even more, we show that $2^\kappa=\kappa^+$ is equivalent to the existence of certain strongly independent families on $\kappa$.

Perhaps the most important result of section one is the fact that the existence of a $\diamondsuit^*$-sequence implies the existence of a strongly independent family on $\omega_1$ of cardinality $2^{\omega_1}$. In this section we also show a relationship between the existence of some of these families and the existence of a strongly inaccessible cardinal.

In the second section, we study the second generalization of independent families, what we have called $\cF$-independent or $\cJ$-independent families, depending on whether $\cF$ is a filter or $\cJ$ is an ideal on a given cardinal $\kappa$. We say that a family is $\cF$-independent (or $\cJ$-independent) if every finite Boolean combination is in $\cF^+$ (or in $\cJ^+$ respectively). For a filter $\cF$ some conditions on it are shown so that there are $\cF$-independent families; in this same direction we show that strongly $\cF$-independent families can also exist, i.e. a kind of double generalization of classical independent families. Later we will focus on the \emph{club} filter, closed and unbounded sets, and show some similarities between this new notion of independence and the classical one. Finally, for an ideal $\cJ\subseteq\cP(\kappa)$, we show that exists a relationship between the existence (or non-existence) of $\cJ$-independent families and the saturation of $\cJ$, therefore with some properties of the cardinal $\kappa$.

\section{Strongly independent families}

For a cardinal $\kappa$ and  $A\subseteq\kappa$, we will use the usual notation, introduced by Shelah in \cite{shelah1992}, $A^0$ denotes $A$ and $A^1$ denotes $\kappa\setminus A$. If $X$ and $Y$ are sets and $s$ is a function, we will use the notation $s;X\rightarrow Y$ to express that $s$ is a partial function from $X$ to $Y$, i.e. $dom(s)\subseteq X$ and $s$ takes its values in $Y$. For a family $\cI$ we will denote the set $\{s;\cI\rightarrow 2: |s|<\omega\}$ by $FF(\cI)$. The rest of the terminology is canonical and it is the one followed by modern literature in set theory.

\begin{mydef}
If $\cI$ is a family of subsets of a cardinal $\kappa$ and $h;\cI\rightarrow 2$, then $\cI^h=\bigcap_{I\in dom(h)}I^{h(I)}$ is \emph{the Boolean combination of $\cI$ determined by} $h$. If $h$ is finite then we say that $\cI^h$ is a finite Boolean combination. If $h$ has cardinality $\lambda$ we say that $\cI^h$ is a \emph{Boolean combination of length} $\lambda$.
\end{mydef}

The set whose elements are all finite Boolean combinations from $\cI$ is the \emph{envelope} of $\cI$ and we denote it by $\ENV(\cI)$.

\begin{mydef}
A family $\cI$ of subsets of a cardinal $\kappa$ is \emph{independent} if every finite Boolean combination of $\cI $ has cardinality $\kappa$.
\end{mydef}

We may generalize independent families allowing larger Boo\-lean combinations.

\begin{mydef}\label{2-3}
A family $\cI$ of subsets of a cardinal $\kappa$ is \emph{strongly independent} if every Boolean combination of length less than $\kappa$ of elements of $\cI$ has size $\kappa$.
\end{mydef}

Normally, after definitions, examples come; instead we now present a typical example of the classical case of an independent family on $\omega$. Latter we shall use it to give examples of the generalizations just introduced.

\begin{example}\label{14}
Let $p_n$ be the $n$-th prime number and $C_n=\{ mp_n : m\in\omega\}$. The family $\cI=\{ C_n : n\in\omega \}$ is independent. 
\end{example}

The family in the previous example is an independent family such that $\cI^h=\emptyset$ for any infinite Boolean combination $h;\omega\rightarrow2$ such that $h^{-1}[\{0\}]$ is infinite. Nevertheless, this does not mean that this family is not strongly independent, since in the case of $\kappa=\omega$, independence and strongly independence agree (it also is the unique cardinal where they do). It is easy to observe that for any independent family $\cI$ infinite on $\omega$ there exists $h;\cI\rightarrow2$ infinite such that $\cI^h=\emptyset$.  In general, in Definition \ref{2-3} we restrict ourselves to Boolean combinations of length less than $\kappa$ because if $\cI$ is an independent family of cardinality $\kappa$ on $\kappa$, there is $h:\cI\rightarrow2 $, with $|h|=\kappa$, such that ${\cI}^h=\emptyset $.

The question naturally arises about for which cardinals it exists (or may exist) a strongly independent family and for which cardinals $\kappa$ there exist \emph{large} strongly independent families, that is, of cardinality $2^\kappa$? Fischer and Montoya in \cite{fischer2019higher} gave a partial answer to this question, which has inspired us to use a guessing principle to construct strongly independent families. 

\begin{mydef}\cite{jensen1972fine}
Let $\kappa$ be a regular cardinal. We say that a sequence  $\langle S_\alpha:{\alpha\in\kappa}\rangle$ is a $\diamondsuit^*(\kappa)$-\emph{sequence} if:
\begin{enumerate}[(1)]
    \item For every $\alpha\in\kappa$, we have that $S_\alpha\subseteq \cP(\alpha)$ and $|S_\alpha|<\kappa$.
    \item For every $X\subseteq\kappa$, the set $\{\alpha\in\kappa: X\cap\alpha\in S_\alpha\}$ is club in $\kappa$. 
\end{enumerate}
The existence of a $\diamondsuit^*(\kappa)$-sequence will be denoted simply as $\diamondsuit^*(\kappa)$. 
\end{mydef}

A very well known and easy consequence of $\diamondsuit^*$ is presented in the next proposition.

\begin{proposition}\label{277}
Let $\kappa$ and $\lambda$ be cardinals such that $\lambda<\kappa$ and $\kappa$ is regular. Then $\diamondsuit^*(\kappa)$ implies $2^\lambda\leq\kappa$.
\end{proposition}

This can be used to show the possibility of having many strongly independent families.

\begin{theorem}\label{EMIP}
Let $\kappa$ be an uncountable regular cardinal. Then $\diamondsuit^*(\kappa)$ implies the existence of a strongly independent family on $\kappa$ of cardinality $2^{\kappa}$.
\end{theorem}
\begin{proof}
Let $\langle S_\alpha:{\alpha\in\kappa}\rangle$ be a $\diamondsuit^*(\kappa)$ sequence and let $C$ be defined as follows:
\[
C=\{\langle\gamma,A\rangle:\gamma\in\kappa\wedge A\subseteq S_\gamma\}.
\]
Since $|S_\alpha|<\kappa$ for every $\alpha\in\kappa$, by Proposition \ref{277}, 

\[
|C|=\sum_{\alpha\in\kappa}2^{|S_\alpha|}\leq\sum_{\alpha\in\kappa}\kappa=\kappa
\]
and it is also clear that $\kappa\leq|C|$, we conclude that $|C|=\kappa$. Thus, constructing a strongly independent family on $\kappa$ is equivalent to doing it on $C$.

For every $X\subseteq\kappa$ let $Y_X$ be defined as follows:
\[
Y_X=\{(\gamma,A)\in C: X\cap\gamma\in A\}.
\]
Aiming to prove that $\cI=\{Y_X: X\subseteq\kappa\}$ is strongly independent, set $\{X_i: i\in I_0\},\{Z_j: j\in I_1\}\subseteq \cP(\kappa)$ two disjoint collections, with $|I_0|,|I_1|<\kappa$.

For every pair $i,i'\in I_0$ with $i\neq i'$ let $\gamma_{i,i'}\in\kappa$ be such that 
\[
X_i\cap\gamma_{i,i'}\neq X_{i'}\cap\gamma_{i,i'}.
\]
Observe that if $\gamma\geq \gamma_{i,i'}$ then $X_i\cap\gamma\neq X_{i'}\cap\gamma$; analogously for $j,j'\in I_1$, with $j\neq j'$ let $\alpha_{j,j'}$ be such that 
\[
Z_j\cap\alpha_{j,j'}\neq Z_{j'}\cap\alpha_{j,j'}.
\]
Finally if $i\in I_0$ and $j\in I_1$, let $\beta_{i,j}\in\kappa$ be such that 
\[
X_i\cap\beta_{i,j}\neq Z_j\cap\beta_{i,j}.
\]

If we define $B\subseteq\kappa$ as: 
\[
B=\left\{\gamma_{i,i'}: i,i'\in I_0\wedge i\neq i'\right\}\cup \{\gamma_{j,j'}: j,j'\in I_1\wedge j\neq j'\}\cup \{\gamma_{i,j}: i\in I_0\wedge j\in I_1\}, 
\]
it is clear that $|B|<\kappa$ and, as $\kappa$ is regular, there exists $\gamma_0\in\kappa$ such that $B$ is bounded by $\gamma_0$. Now, if $\gamma\in\kappa$ is larger that $\gamma_0$, then this one satisfies the following:
\begin{enumerate}[(1)]
    \item $X_i\cap\gamma\neq X_{i'}\cap\gamma$ if $i,i'\in I_0$ with $i\neq i'$.
    \item $Z_j\cap\gamma\neq Z_{j'}\cap\gamma$ if $j,j'\in I_1$ with $j\neq j'$.
    \item $X_i\cap\gamma\neq Z_{j}\cap\gamma$ if $i\in I_0$ with $j\in I_1$.
\end{enumerate}

For every $i\in I_0$, consider $D_i=\{\gamma\in\kappa: X_i\cap\gamma\in S_\gamma\}$, which is a club, now put $D=\bigcap_{i\in I_0}D_n$ and let $\gamma\in D$ such that $\gamma>\gamma_0$.

Let $A_{\gamma}\subseteq S_\gamma$ be defined as:
\[
A_\gamma=\{X_i\cap\gamma: i\in I_0\}.
\]
So we have that $(\gamma,A_\gamma)\in Y_{X_i}$ for every $i\in I_0$ and $(\gamma, A_\gamma)\not\in Y_{Z_j}$ for every $j\in I_1$. This proves that:
\[
(\gamma,A_\gamma)\in \bigcap_{i\in I_0}Y_{X_i}\setminus\bigcup_{j\in I_1}Y_{Z_j}
\]
and as this happens for every $\gamma\in D$ such that $\gamma>\gamma_0$, then:
\[
\big|\bigcap_{i\in I_0}Y_{X_i}\setminus\bigcup_{j\in I_1}Y_{Z_j}\big|=\kappa,
\]
which finishes the proof.
\end{proof}

If $\kappa $ is strongly inaccessible then $\langle\mathscr P(\alpha):{\alpha\in\kappa}\rangle$ turns out to be a $\diamondsuit^*(\kappa)$-sequence, hence the previous theorem in particular implies that for every strongly inaccessible cardinal there is a {large} strongly independent family on it, which is a result obtained by Fischer and Montoya in \cite{fischer2019higher}; however, Theorem \ref{EMIP} gives a broader spectrum of cardinals for which there are consistently strong independent {large} families on them. For example, under $\mathbf{V=L}$, Theorem \ref{EMIP}, implies that large strongly independent families exist on many cardinals. In fact in \cite {jensen1972fine} Jensen proved:

\begin{theorem}\cite{jensen1972fine}
$\mathbf{V=L}$ implies $\diamondsuit^*(\kappa)$ for every successor cardinal $\kappa$. 
\end{theorem}
\begin{corollary}
$\mathbf{V=L}$ implies that for every successor cardinal $\kappa$ it exists a strongly independent family of cardinality $2^\kappa$.
\end{corollary}

On the other hand, the existence of strongly independent families on successor cardinals is also closely related to the Generalized Continuum Hypothesis. 

\begin{theorem}\label{19}
Let $\kappa$ be an infinite cardinal. The following two conditions are equi\-valent.
\begin{enumerate}[(1)]
    \item There is a strongly independent family on $\kappa^+$ of cardinality $\kappa$.
    \item The equality $2^\kappa=\kappa^+$ is true.
\end{enumerate}
\end{theorem}
\begin{proof}
$(1)\Rightarrow (2)$. Let $\cI=\{X_\alpha:\alpha\in\kappa\}$. Now, since $\cI$ is a strongly independent family on $\kappa^+$, for all $h\in2^\kappa $ we have that $\cI^h$ has cardinality $\kappa^+$ and it is clear that if $h, g\in2^\kappa$ are different then $\cI^h$ and $\cI^g$ are disjoint. For every $h\in 2^\kappa$, let $x_h\in\cI^h$; then the set $\{x_h: h\in2^\kappa\}$ is a subset of $\kappa^+$ and has cardinality $2^\kappa$, so $2^\kappa\leq\kappa^+$ and therefore $2^\kappa=\kappa^+$.

$(2)\Rightarrow (1)$. Let $f:\kappa^+\rightarrow 2^\kappa\times \kappa^+$ be a bijection (considering $2^\kappa$ as the set of all functions from $\kappa$ to $2$). For every $h\in 2^\kappa$, let $X_h=f^{-1}(\{h\}\times\kappa^+)$ and for every $\alpha\in\kappa$ let $I_\alpha$ be defined as follows:
\[
I_\alpha=\bigcup\left\{X_h: h\in 2^\kappa\mm\{\overline{1}\}\y h(\alpha)=0\right\},
\]
where $\overline{1}$ denotes the function $f:\kappa\rightarrow 2$ with constant value $1$.

Let $\cI=\{I_\alpha:\alpha\in\kappa\}$. It is clear that if $h\in 2^\kappa\mm\{\overline{1}\}$ then $\cI^h\supseteq X_h$ and, as $|X_h|=\kappa^+$, we have that $|\cI^h|=\kappa^+$, which proves that $\cI$ is strongly independent. 
\end{proof}

The following results are simple corollaries of Theorem \ref{19}.

\begin{corollary}
There exists an infinite strongly independent family on $\omega_1$ if and only if $\mathsf{CH}$ is satisfied, thus, the existence of an infinite strongly independent family on $\omega_1$ is independent from $\mathsf{ZFC}$.
\end{corollary}

\begin{corollary}
Let $\kappa$ be an inaccessible cardinal (limit and regular) such that for e\-ve\-ry infinite cardinal $\lambda<\kappa$ it exists a strongly independent family on $\lambda^+$ of cardinality $\lambda$, then $\kappa$ is strongly inaccessible.
\end{corollary}

\begin{proof}
We only need to verify that $\kappa$ is a strong limit cardinal. Let $\lambda\in\kappa$; as $\kappa$ is limit it follows that $\lambda^+<\kappa$. On the other hand, as it exists a strongly independent family of size $\lambda$ on $\lambda^+$, then $2^\lambda=\lambda^+$ and so $2^\lambda<\kappa$, which finishes the proof.
\end{proof}

\begin{corollary}
If $\kappa$ is inaccessible and for every $\lambda<\kappa$ there is a strongly independent family of cardinality $\lambda$ on $\lambda $, then $\kappa$ is strongly inaccessible.
\end{corollary}

Although we already know some sufficient conditions for the existence of strongly independent families, an interesting property of these is that they do not satisfy the conditions to apply Zorn's Lemma (unlike the classical independent families), which is the standard way to prove that maximal objects with some property exist. It is therefore of great interest to know:

\begin{preg}\label{Q}
For which cardinals are there strongly independent maximal families on them?
\end{preg}

It is not known yet if these families exist for any cardinal, the only results we have so far are in the direction of the \textit{not maximality}.

\begin{mydef}
A strongly independent family $\cI$ on a cardinal $\kappa$ is \emph{maximal} if there is no other strongly independent family on $\kappa$ that properly extends it.
\end{mydef}

\begin{theorem}
On any infinite cardinal $\kappa>\omega$ there exists an independent family that is not strongly independent. 
\end{theorem}

\begin{proof}
We know that there exists a bijection between $\kappa$ and $\omega\times\kappa$, so we are going to construct the desired independent family on $\kappa\times\omega$. For every $n\in\omega$ let $I_n=\kappa\times C_n$, where the $C_n$ are as in the Example \ref{14}, and let $\cI=\{I_n: n\in\omega\}$.

Clearly if $h;\omega\rightarrow 2$ is finite, then for every $\alpha\in\kappa$ we have that $(\{\alpha\}\times\omega)\cap \cI^h$ is infinite, in particular $\cI^h$ has size $\kappa$. On the other hand, if $h:\omega\rightarrow 2$ is such that $h^{-1}[\{0\}]$ is infinite, then for every $\alpha\in\kappa$ we have that $(\{\alpha\}\times\omega)\cap\cI^h=\emptyset$, which implies that $\cI^h=\emptyset$, thus $\cI$ is not strongly independent.
\end{proof}

Note that the family constructed in the proof of the previous theorem can be extended to a maximal independent family $\cJ$, and since $\cI\subseteq\cJ$, then $\cJ$ is not strongly independent either. Thus we have the next corollary. 

\begin{corollary}
For every infinite cardinal $\kappa$ there exists a maximal independent family on $\kappa$ that is not strongly independent.
\end{corollary}

As in the classical case of independent families, we know that strongly independent families small in cardinality are not maximal.

\begin{proposition}\label{111}
If $\cI$ is a strongly independent family on a cardinal $\kappa$ such that $|\cI|<\kappa$, then $\cI$ is not maximal.
\end{proposition}

\begin{proof}
Let $\cI=\{I_\alpha: \alpha\in\lambda\}$ with $\lambda<\kappa$ and for each $h:\lambda\rightarrow2$ let $X_h=\cI^h$. Now each set $X_h$ is of cardinality $\kappa$ and if $ h, g \in 2^\lambda $ are different then $X_h\cap X_g=\emptyset $, this implies that $ 2^\lambda\leq\kappa $. Let $\langle Y_\alpha:{\alpha\in\kappa}\rangle$ be an enumeration of $\{X_h: h \in2^\lambda\}$ such that every $X_h$ appears $\kappa$ times. Let $ a_0, b_0 \in Y_0 $ be such that $ a_0<b_0 $ and suppose that $ a_\beta $ and $b_\beta $ have been already defined for all $ \beta<\alpha $. Since $Y_\alpha$ has cardinality $\kappa$ there are $ a_\alpha, b_\alpha \in Y_\alpha $ such that for all $ \beta\in\alpha $ it holds that $ a_\beta, b_\beta<a_\alpha $ and also $a_\alpha <b_\alpha $. Now let $Z=\{a_\alpha : \alpha \in \kappa \} $. By the construction of $ Z $ we have that $ Z\cap X_h $ and $ (\kappa\setminus Z) \cap X_h $ have cardinality $\kappa $ for all $ h \in2^\lambda $, that is, $ \cI\cup\{Z\}$ is a strongly independent family.
\end{proof}

Note that the above proof is not applicable to strongly independent families of cardinality $\kappa$.

The following shows, in the same direction of Proposition \ref{111},  that another class of strongly independent families are not maximal neither.

\begin{mydef}
Let $\kappa$ be an infinite cardinal.
\begin{enumerate}[(1)] 
    \item Let $\cF\subseteq\cP(\kappa)$ and $X\subseteq\kappa$, we say that $X$ \emph{splits} $\cF$ if $Y\cap X$ and $Y\setminus X$ have size $\kappa$ for all $Y\in\cF$. 
    \item A family $\cR\subseteq\cP(\kappa)$ is \emph{unsplittable} (or \emph{reaping}) if there is not $X\subseteq \kappa$ that splits $\cR$. 
    \item $\mathfrak{r}(\kappa)$ is the smallest cardinality of a unsplittable family on $\kappa$.
\end{enumerate}
\end{mydef}

\begin{theorem}
(Fischer-Montoya \cite{fischer2019higher}) Let $\kappa$ be an infinite regular cardinal. If $\cI$ is a strongly independent family on $\kappa$ such that $|\{\cI^h: h;\cI\rightarrow 2 \wedge |h|<\kappa\}|<\mathfrak{r}(\kappa)$ then $\cI$ is not maximal.
\end{theorem}

In \cite{kunen}, K. Kunen studied maximal $\sigma$-independent families in uncountable cardinals; that is, maximal families which are independent with respect to Boolean combinations of countable length. For instance he proved 

\begin{theorem}[Kunen]
If $\mathsf{ZFC}$ plus the existence of a measurable cardinal is consistent, so is $\mathsf{ZFC}$ plus the existence of a maximal $\sigma$-independent family $\mathscr{S}\subset\mathscr{P}(\omega_1)$.
\end{theorem}

His methods are ad hoc and it does not seem possible to generalized them to answer Question \ref{Q}; however, this gives an idea of the consistency strength one has to face to answer Question \ref{Q}. Kunen's paper also shows all the complexity of the property of maximality for independent families on uncountable cardinals. As we said earlier, we were unable to present properties that guarantee maximality for strong independent families.  In the next section we take a different approach to generalized the classical case. Again the property of being maximal for those is perhaps even harder. For example, we were unable to prove that a countable $\cC$-independent family cannot be maximal.  See Theorem  \ref{Cnumerable}.

\section{$\cF$-independent families}

Let $\cF$ be a filter on $\kappa$. A subset $X\subseteq\kappa$ is $\cF$-positive if $X\cap Y\neq\emptyset$ for every $Y\in\cF$; we denote the family of $\cF$-positive subsets by $\cF^+$. If $\cJ\subseteq\cP(\kappa)$ is an ideal then $\cJ^+=\{X\subseteq\kappa: X\not\in\cJ\}$.

If $\cF$ a filter on a cardinal $\kappa$, we denote by $\cF^*$ its dual ideal, i.e, the ideal $\{X\subseteq \kappa:\kappa\setminus X\in\cF\}$.

\begin{mydef}
A family $\cI\subseteq\cP(\kappa) $ is $\cF$-\emph{independent} if every finite Boolean combination of $\cI$ is in $\cF^+$. Similarly if $\cJ$ is an ideal then $\cI$ is $\cJ$-\emph{independent} if every finite Boolean combination of $\cI$ is in $\cJ^+$.
\end{mydef}

Note that a family is $\cF$-independent if and only if it is $\cF^*$-independent. On the other hand, if $\mathcal{F}_r$ is the Fréchet filter, then a family is $\cF_r$-independent if and only if it is independent. It is also clear that if $ \cI $ is a $\cF$-independent family on $\kappa$ and $X\in\cI $, then $X$ is \emph{$\cF$-double positive}, that is, $X, \kappa\setminus X\in\cF^+$, consequently if $\cF $ is an ultrafilter, there are no $\cF$-independent families. The natural question is to know for which filters (or ideals) $\cF\subseteq\cP(\kappa)$ (in addition to the Fréchet's one) there is a $\cF$-independent family.

\begin{proposition}
Let $\cF$ be a filter of the form $\cF=\{A\subseteq\kappa: B\subseteq A\}$ for some $B\subseteq\kappa$.
\begin{enumerate}[(1)]
    \item If $B$ is finite then there are not $\cF$-independent infinite families, furthermore, if $|B|=n$ then there are not $\cF$-independent families of cardinality $n$. 
    \item If $|B|=\lambda$ with $\lambda$ infinite, there exist $\cF$-independent families of cardinality $2^\lambda$ but not of cardinality $(2^\lambda)^+$.
\end{enumerate}
\end{proposition}

\begin{proof}
$(1)$ Note that $\cF^+=\{X\subseteq\kappa: X\cap B\neq\emptyset\}$. Let $B=\{x_0,\dots,x_{n-1}\}$ and suppose that $X_0,\dots,X_{n-1}\in \cI$ are all distinct, where $\cI$ is an $\cF$-independent family. For each $i\in n$, if $x_i\in X_i$ let $h(i)=1$ and $h(i)=0$ otherwise; so we have that $x_i\not\in X_i^{h(i)}$. 
Then for every $x\in B$ we have that:
\[
x\not\in \bigcap_{i\in n} X_i^{h(i)}=\cI^h,
\]
so $\cI^h\cap B=\emptyset$ and therefore $\cI^h\not\in\cF^+$, which contradicts the fact that $\cI$ is $\cF$-independent.

$(2)$ Again note that $\cF^+=\{X\subseteq\kappa: X\cap B\neq\emptyset\}$. Now let $\cI=\{X_\alpha: \alpha\in 2^\lambda\}$ be an independent family of subsets of $B$ and for each $\alpha\in 2^\lambda$ let $Y_\alpha=X_\alpha\cup (\kappa\setminus B)$ and let $\widehat{\cI}=\{Y_\alpha:\alpha\in 2^\lambda\}$. Clearly if $h;2^\lambda\rightarrow 2$ is finite then $\cI^h\subseteq\widehat{\cI}^h$ and as $\cI$ is independent on $B$ we have that:
\[
\emptyset\neq B\cap\cI^h=B\cap\widehat{\cI}^h,
\]
which proves that $\widehat{\cI}^h\in\cF^+$, therefore $\widehat{\cI}$ is $\cF$-independent.

If $\cI\subseteq\kappa$ has cardinality $(2^\lambda)^+$, as $|B|=\lambda$, there exist $X,Y\in\cI$ distinct such that $X\cap B=Y\cap B$, but then $(X\setminus Y)\cap B=\emptyset$, which proves that $X\setminus Y\not\in\cF^+$, thus $\cI$ is not $\cF$-independent.
\end{proof}

As anticipated, the two generalizations of independence studied in this work are compatible with each other, that is, we can \emph{merge} the two notions in order to obtain families with more combinatorial properties.

\begin{mydef}
Let $\cF\subseteq\cP(\kappa)$ be a filter (respectively $\cJ\subseteq\cP(\kappa) $ an ideal). A family $\cI\subseteq \cP(\kappa)$ is \emph{strongly ${\cF}$-independent} (respectively \emph{strongly $\cJ$-independent}) if every Boolean combination of length less than $\kappa $ of $\cI$ is in $\cF^+$ (respectively in $\cJ^+$).
\end{mydef}

We will study a little more of these families below.

\subsection{$\cC$-independent families}

For each regular cardinal $\kappa$ let $\cC_{\kappa}\subseteq\cP(\kappa)$ be the club filter, that is, the filter generated by closed and unbounded sets (when the context is clear we will call $\cC_\kappa$ simply as $\cC$). $\cC_{\omega_1}$ is a very important filter in the study of the combinatorics of $\omega_1$, therefore a couple of questions  arise naturally: Are there $\cC_{\omega_1}$-independent families?  Is every maximal $\cC$-independent family strongly $\cC$-independent? Answers to these questions can be found in Proposition \ref{212} and Corollary \ref{213}, respectively.

First of all, let us note that as for every filter $\cF$, the union of $\cF$-independent families is an $\cF$-independent family, then if there are $\cF$-independent families then there are maximal ones (by Zorn's Lemma).

Remember that $\cC$-positive sets are called \emph{stationary} sets; one of the most important results about stationary sets is the following:

\begin{lemma}[\cite{solovay1971real}, \cite{kanamori2008higher}]
For each uncountable regular cardinal $\kappa$ we have that $\kappa$ is the union of as many as $\kappa$ disjoint stationary sets.
\end{lemma}

\begin{corollary}\label{3.8}
For each uncountable regular cardinal $\kappa$ and each $\lambda\leq\kappa$ we have that $\kappa$ is the union of $\lambda $ disjoint stationary sets.
\end{corollary}

The following two results are consequences of this last corollary; their proof follow the scheme of the proof of Theorem \ref{19}.

\begin{proposition}\label{212}
For every uncountable regular cardinal $\kappa$ there exists an infinite $\cC$-independent family.
\end{proposition}

\begin{proof}
By Corollary \ref{3.8} there is a countable collection $\{X_s: s\in2^{<\omega}\}$ of disjoint stationary subsets whose union is $\kappa$, say indexed with the set $2^{<\omega}$.

Now, for every $n\in\omega$, let $I_n\subseteq\kappa$ be defined as follows:
\[
I_n=\bigcup\{X_s: s\in2^{<\omega}\wedge n\in dom(s)\wedge s(n)=0\}.
\]

It turns out that $\cI=\{I_n: n\in\omega\}$ is a  $\cC$-independent family, since every finite Boolean combination of $\cI$ contains some combination of the form
\[
\bigcap \{I_n^{s(n)}:{n\in dom(s)}\}
\]
for some $s\in 2^{<\omega}$ and also:
\[
X_s\subseteq\bigcap\{I_n^{s(n)}:{n\in dom(s)}\},
\]
which proves that every finite Boolean combination of $\cI$ contains a stationary set, therefore is stationary.
\end{proof}

\begin{corollary}\label{213}
For any cardinal $\kappa\geq\omega_1$ it exists a $\cC$-independent maximal family on $\kappa$ that is not strongly $\cC$-independent.
\end{corollary}

\begin{proof}
Let $\{X_m: m\in\omega\}$ be a partition of $\kappa$ into stationary sets. Now for every $n\in\omega$ let $Y_n=\bigcup\{X_m: m\in C_n\}$, where the $C_n$ are as in the Example \ref{14}. Consider $\cI=\{Y_n: n\in\omega\}$; then it is easy that $\cI$ is $\cC$-independent but for $h;\omega\rightarrow2$ such that $h^{-1}[\{0\}]$ is infinite we have that $\cI^h=\emptyset$, which proves that $\cI$ is not strongly $\cC$-independent. Extending $\cI$ to a maximal $\cC$-independent family the result is obtained. 
\end{proof}

\begin{theorem}\label{310}
The following statements are equivalent for a cardinal $\kappa$:
\begin{enumerate}
    \item $2^\kappa=\kappa^+$.
    \item There exists a strongly independent family on $\kappa^+$ of size $\kappa$.
    \item There exists a strongly $\cC$-independent family on $\kappa^+$ of size at least $\kappa$.
\end{enumerate}
    \end{theorem}

\begin{proof}
We only prove (1)$\Rightarrow$ (3). Let $\{X_f: f\in2^\kappa\}$ be a partition of $\kappa^+$ into stationary sets and for every $\alpha\in\kappa$ let $I_\alpha$ be defined by
\[
I_\alpha=\bigcup\left\{X_f: f\in2^\omega\y(f(\alpha)=0)\right\}.
\]
Let $\cI=\{I_\alpha: \alpha\in\kappa\}$. It is clear that if $f;\kappa\rightarrow2$ then $\cI^f\supseteq X_h$ for some $h\in2^\kappa$ and, as $X_h$ is stationary, $\cI^f$ is stationary too, which proves that $\cI$ is strongly $\cC$-independent.
\end{proof}

We now know that there are countable $\cC$-independent families on $\omega_1$. Are there uncountable $\cC$-independent families on $\omega_1$? Furthermore, are there $\cC$-independent families of cardinality $2^{\omega_1}$?

\begin{theorem}
Let $\kappa$ and $\lambda$ be cardinals such that $\omega\leq\lambda\leq2^{\kappa}$ and $\kappa$ is regular. Then, on $\kappa$, there is a $\cC$-independent family of cardinality $\lambda$.
\end{theorem}

\begin{proof}
Let $\{X_\beta: \beta\in\kappa\}$ be a partition of $\kappa$ into stationary sets. Now let $\cI=\{I_\alpha : \alpha\in\lambda\}$ be an independent family of cardinality $\lambda $ on $\kappa $. For every $\alpha\in\lambda$, let $\widehat{I_\alpha}\subseteq\kappa$ be defined as follows:
\[
\widehat{I_\alpha}=\bigcup\{X_\beta:\beta\in I_\alpha\}.
\]
Now let $\widehat{\cI}=\{\widehat{I_\alpha}:\alpha\in\lambda\}$. Clearly $\widehat{\cI}$ has size $\lambda$, then the only thing left to prove is that it is a $\cC$-independent family. Let $s;\lambda\rightarrow2$ be finite, we want to see that $\widehat{\cI}^s $ is stationary. Since $\cI$ is independent there is $\beta\in\cI^s$, but this means that if $s(\alpha)=0$ then $X_\beta\subseteq\widehat{I_\alpha}$ and if $s(\alpha)=1$ then $X_\beta\cap\widehat{I_\alpha}=\emptyset$, that is, $X_\beta\subseteq\widehat{\cI}^s$, and since $X_\beta$ is stationary $\widehat{\cI}^s$ is also stationary.
\end{proof}

As in the classical case of independent families, one would expect that the coun\-table $\cC$-independent families are not maximal; however, it seems complicated to establish that. Our ideas about generalizing the classical proof, doing a disjoint refinement of the envelope or using a $\diamondsuit^\sharp$-sequence have failed. The following is a modification of the main construction from \cite{HS}.

\begin{theorem}\label{Cnumerable}
Under $\mathbf{V=L}$, a countable $\cC_{\omega_1}$-independent family is not maximal.
\end{theorem}

\begin{proof}
Let $\cI$ be a countable $\cC$-independent family, and let $\{E_n:\nom\}$ be an enumeration of its envelope. For each limit ordinal $\gamma<\omega_1$ set
\[
\mathscr A_\gamma=\{\alpha<\omega_1 : L(\alpha)\models\mathsf{ZF^-}\y\gamma=\omega_1^{L(\alpha)}\}.
\]
Since $\{\varrho<\omega_1 : L(\varrho)\prec L(\omega_1)\}$ is unbounded in $\omega_1$, it follows that $\mathscr A_\gamma$ is at most countable for each limit $\gamma<\omega_1$. It is also known that $\{\gamma<\omega_1 : \mathscr A_\gamma\neq\emptyset\}$ contains a club.   Let
\[
\mathscr G_\gamma=\{C\subseteq\gamma : C\mbox{ is club in }\gamma\y(\exists\alpha\in\mathscr A_\gamma)(C\in L(\alpha))\}.
\]
Then $\mathscr G_\gamma$ is countable and since $\mathsf{ZF}^-$ suffices to prove that the intersection of a finite collection of club subsets is a club subset, it follows that $\mathscr G_\gamma$ is closed under finite intersections.

Consider as well 
\[
\mathscr S_\gamma=\{S\subseteq\gamma : (\exists\alpha\in\mathscr A_\gamma)(S\in L(\alpha))\y(\forall C\in\mathscr G_\gamma)(C\cap S\neq\emptyset) \}.
\]
Once again $\mathscr S_\gamma$ is countable; fix an enumeration $\{S_n:\nom\}$ of $\mathscr S_\gamma$ in which each element appears infinitely often and some simple enumeration $\{C_n:\nom\}$ of $\mathscr G_\gamma$.  Now consider a cofinal sequence $\langle\alpha_n:\nom\rangle$ in $\mathscr A_\gamma$ such that 
\[
S_n\in L(\alpha_n)\y(\forall m\leq n)(C_m\in L(\alpha_n)).
\]

Since $L(\alpha_0)\models\mbox{``}S_0\mbox{ is stationary in }\gamma\mbox{''}$ pick
\[
\xi_0\in S_0\cap C_0\quad\mbox{and}\quad\eta_0\in S_0\cap C_0\mm(\xi_0+1),
\]
and recursively 
\[
\xi_{n+1}\in(S_{n+1}\cap\bigcap_{k\leq n+1} C_{k})\mm(\eta_n+1)\quad\mbox{and}\quad\eta_{n+1}\in(S_{n+1}\cap\bigcap_{k\leq n}C_k)\mm(\xi_{n+1}+1),
\]
for all $\nom$. This way we have built two disjoint subsets $A_\gamma=\{\xi_n:\nom\}$ and $B_\gamma=\{\eta_n:\nom\}$.

Put $A=\bigcup\{A_\gamma:\gamma\in\mathrm{Lim}(\omega_1)\}$ and $B=\bigcup\{B_\gamma:\gamma\in\mathrm{Lim}(\omega_1)\}$.\bigskip

\noindent\emph{Claim}: $(\forall k\in\omega)(E_k\cap A\mbox{ is stationary})$.\medskip

Fix a club subset $C\subseteq\omega_1$. Define recursively a sequence of elementary submodels $M_\nu\prec L(\omega_2)$ for $\nu<\omega_2$ as follows:
\begin{itemize}
    \item $M_0$ is the smallest $M\prec L(\omega_2)$ such that $\{E_n:\nom\},C\in M$,
    \item $M_{\nu+1}$ is the smallest $M\prec L(\omega_2)$ such that $M_\nu\cup\{M_\nu\}\subseteq M$,
    \item $M_\xi=\bigcup_{\nu<\xi}M_\nu$ whenever $\xi$ is a limit ordinal.
\end{itemize}

By the Condensation Lemma,  $M_\nu\cap L(\omega_1)$ is transitive, set $\alpha_\nu=M_\nu\cap\omega_1$.  Then $\langle\alpha_\nu:\nu<\omega_1\rangle$ 
is a normal sequence in $\omega_1$. Use Mostowski's Collapse $\pi_\nu:M_\nu\cong L(\beta_\nu)$ to get
\begin{itemize}
    \item $\pi_\nu\rtt L(\alpha_\nu)=\id\rtt L(\alpha_\nu)$,
    \item $\pi_\nu(\omega_1)=\alpha_\nu$,
    \item $\pi_\nu(C)=C\cap\alpha_\nu$,
    \item $(\forall\nom)(\pi_\nu(E_n)=E_n\cap\alpha_\nu)$.
\end{itemize}
Consider the set $K$ of limit points of $\langle\alpha_\nu:\nu<\omega_1\rangle$. Obviously $K$ is a club in $\omega_1$, if $\gamma\in K$, then 
\[
\gamma=\sup_{\nu<\zeta}\alpha_\nu=\sup_{\nu<\zeta}\beta_\nu,
\]
for some ordinal $\zeta<\omega_1$, and hence $\gamma=\alpha_\zeta$. To see this, it is enough to show $\alpha_\nu<\beta_\nu<\alpha_{\nu+1}$.  Clearly $\alpha_\nu<\beta_\nu$.  Since $\beta_\nu$ is definable from $M_\nu$ as $L(\beta_\nu)$ is the transitive collapse of $M_\nu$ and that definition relativises to $L(\omega_2)$. Thus $\beta_\nu\in M_{\nu+1}$ as $M_\nu\in M_{\nu+1}\prec L(\omega_2)$.  Henceforth $\beta_\nu\in\alpha_{\nu+1}$.

Note that $\beta_\zeta\in\mathscr A_\gamma$ since $L(\beta_\zeta)\models\gamma=\omega_1$ and $L(\beta_\zeta)\models\mathsf{ZF}^-$.  Thus $C\cap\gamma=\pi_\zeta(C)\in L(\beta_\zeta)$ and of course $L(\beta_\zeta)$ models that $\pi_\zeta(C)$ is a club in $\gamma$.  This implies $C\cap\gamma\in\mathscr G_\gamma$. It is also true that $E_k\cap\gamma\in L(\beta_\zeta)$, then $E_k\cap\gamma=S_n\in\mathscr S_\gamma$, for infinitely many $\nom$.  Since $A_\gamma\mm(C\cap\gamma)$ is finite and $A_\gamma$ is built in such a way that $A_\gamma\cap(E_k\cap\gamma)$ is infinite, this shows that $E_k\cap A$ is stationary in $\omega_1$.

Analogously $E_k\cap B$ is stationary in $\omega_1$ for all $k\in\omega$.  It follows that $\cI\cup\{A\}$ is also $\cC$-independent.
\end{proof}

Observe that it is easily possible that $A\cap B\neq\emptyset$; however, it is not hard to show that $A\mm B$ and $B\mm A$ are stationary as well.

Maximal $\cC$-independent families have many properties analogous to those of maximal independent ones in the classical case. For example, it is easy to prove that if $\cI$ is $\cC$-independent and finite then it is not maximal. Indeed, let us say that $\cI=\{I_i: i\in n \}$ for some $n\in\omega $ and note that for each $s: n\rightarrow2$, the set $\cI^s=\bigcap_{i \in n} I^{s(i)} $ is stationary; furthermore, if $ s, t: n\rightarrow2 $ are different, $ \cI^s$ and $\cI^t$ are disjoint. For each $s:n\rightarrow2$ let $A_s$ and $B_s$ be a partition of $\cI^s $ into two disjoint stationary sets and let $A=\bigcup\{A_s : s\in2^n \}$ . It is clear that $ A\not\in\cI $ and $\cI\cup\{A \}$ is $\cC$-independent.  
Note that, since we can always split a stationary set into two stationary subsets, the above guarantees that we can recursively construct $\cC$-independent families of cardinality $n\in\omega$ and thus obtain a countable $\cC$-independent family. The advantage of this method is that it only requires the fact that a stationary set can be split into two stationary sets and not necessarily into infinite ones. 

\subsubsection{Dense $\cC$-independent families}

All the properties shown next for $\cC$-inde\-pen\-dent families were proved for the classic independence by Goldstern and Shelah in \cite{goldstern1990ramsey}, this proves that $\cC$-independent families on $\omega_1$ behave similarly as the independent ones on $\omega$.

\begin{mydef}
If $\cI$ is a $\cC$-independent family then we define the \emph{ideal associated} to $\cI$ as: 
\[
\cJ_\cI=\{A\subseteq\omega_1:(\forall f\in FF(\cI))(\exists g\in FF(\cI))(g\supseteq f \wedge \cI^g\cap A\text{ is not stationary})\}.
\]
Clearly $\cJ_\cI$ is an ideal that contains the ideal of the non-stationary sets.
\end{mydef}

\begin{mydef}
\begin{enumerate}[(1)]
    \item If $X,Y\subseteq\omega_1$, we say that $X$ is \emph{almost contained} in $Y$ if $X\setminus Y$ is not a stationary set and we denote this by $X\subseteq_{\footnotesize\textsc{ns}}Y$.
    \item For a family $\cX$ of subsets of $\omega_1$ and $Y\subseteq\omega_1$, we say that $Y$ is \emph{pseudointersetion} of $\cX$ if $Y\subseteq_{\footnotesize\textsc{ns}}X$ for every $X\in\cX$.
\end{enumerate}
\end{mydef}

\begin{mydef}
A $\cC$-independent maximal family is \emph{dense} if for every $A\in{\cJ_\cI}^+$ it exists $g\in FF(\cI)$ such that $\cI^g\subseteq_{\footnotesize\textsc{ns}}A$.
\end{mydef}

This can be interpreted as follows: a $\cC$-independent family is dense if the envelope of $\cI$ is a \emph{base} of ${\cJ_\cI}^+$; let us also note that for all $f\in FF(\cI)$ we have that $\cI^f\in {\cJ_\cI}^+$, since $f$ itself is a witness of this.  
Next we use the following standard notation, if $\mathscr{A}\subset\mathscr{P}(X)$ and $Y\subset X$, then $\mathscr{A}\restriction Y$ is the family $\{A\cap Y:A\in\mathscr A\}$.

\begin{proposition}
If $\cI$ is a maximal $\cC$-independent family, it exists $f\in FF(\cI)$ such that for every $g\in FF(\cI)$, with $g\supseteq f$, $\cI\restriction\cI^g$ is maximal.
\end{proposition}

\begin{proof}
Let $\{f_n: n\in\omega\}$ be a maximal family with the following properties:
\begin{enumerate}[(1)]
    \item If $n\neq m$, $f_n$ and $f_m$ are incompatible.
    \item $\cI\restriction\cI^{f_n}$ is not maximal for every $n\in\omega$.
\end{enumerate}
Note that by condition 1) and since $FF(\cI)$ is ccc, this collection is at most countable (in principle it could be finite but assume without loss of generality that it is countable).

Now, for every $n\in\omega$ let $A_n\subseteq\cI^{f_n}$ be such that $\cI\restriction\cI^{f_n}\cup\{A_n\}$ is $\cC$-independent on $\cI^{f_n}$ and let $A=\bigcup_{n\in\omega}A_n$. Since $\cI$ is maximal it exists $f\in FF(\cI)$ such that $\cI^f\cap A$ or $\cI^f\setminus A$ is not stationary. Let us suppose without loss of generality that $\cI^f\cap A$ is not stationary. 
We claim that $f$ is incompatible with every $f_n$; to see this, suppose that $f$ and $f_n$ are compatible, i.e. suppose that $f\cup f_n$ is a function. Thus $\cI^{f\cup f_n}\in \ENV(\cI\restriction\cI^{f_n})$, in particular we have that:
\[
\cI^f\cap A\supseteq\cI^{f\cup f_n}\cap A\supseteq\cI^{f\cup f_n} \cap A_n,
\]
but this is impossible, since in that case $\cI^{f\cup f_n} \cap A_n$ is stationary as $\cI^f\cap A$ is not.

Since $f$ is incompatible with every $f_n$ then so is every $g\in FF(\cI)$ such that $g\subseteq f$, therefore $\cI\restriction\cI^g $ is maximal, otherwise the maximality of $\{f_n: n\in\omega\}$ would be contradicted.
\end{proof}

\begin{lemma}
If $\cI$ is a $\cC$-independent maximal family such that for every $f\in FF(\cI)$ the family $\cI\restriction\cI^f$ is maximal, then $\cI$ is dense.
\end{lemma}

\begin{proof}
Let $A\in {\cJ_\cI}^+$, this means that there exists $f\in FF(\cI)$ such that for every $ g\in FF(\cI) $ that extends to $f$ we have that $\cI^g\cap A$ is stationary. As $\cI\restriction\cI^f$ is maximal, it exists $g\in FF(\cI)$, $g\supseteq f$ such that either $\cI^g\cap A$ or $\cI^g\setminus A$ is not stationary, but we know that $\cI^g\cap A$ is stationary, then necessarily $\cI^g\setminus A$ is not, i.e. $\cI^g\subseteq_{\footnotesize\textsc{ns}}A$, which is what we wanted. 
\end{proof}

\begin{proposition}
If $\cI$ is a maximal $\cC$-independent family which is dense, then $\rfrac{P(\omega_1)}{\cJ_\cI}$ is ccc. 
\end{proposition}

\begin{proof}
By contradiction. Suppose that $\{X_\alpha:\alpha\in\omega_1\}\subseteq{\cJ_\cI}^+$ is such that if $\alpha\neq\beta$ then $X_\alpha\cap X_\beta\in\cJ_\cI$. Since $\cI$ is a dense family, for every $\alpha\in\omega_1$ it exists $f_\alpha\in FF(\cI)$ such that $\cI^{f_\alpha}\subseteq_{\footnotesize\textsc{ns}} X_\alpha$. Now if $\alpha\neq\beta$ then $f_\alpha$ and $f_\beta$ are incompatible, otherwise $\cI^{f_\alpha\cup f_\beta}=\cI^{f_\alpha}\cap\cI^{f_\beta} \subseteq_{\footnotesize\textsc{ns}}X_\alpha\cap X_\beta\in\cJ_\cI$. But now $\cI^{f_\alpha\cup f_\beta}\in\cJ_\cI$ (as $\cJ_\cI$ contains the non-stationary sets), and this is a contradiction as $\ENV(\cI)\subseteq {\cJ_\cI}^+$. 

Thus the family $\{f_\alpha:\alpha\in\omega_1\}$ is an antichain in $FF(\cI)$, but this contradicts the fact that $FF(\cI)$ is ccc.
\end{proof}

\subsubsection{Strongly $\cC$-independent families}

\begin{lemma}
Let $\mathcal{E}=\{E_n: n\in\omega\}$ be a nested collection of stationary sets, i.e. $E_{n+1}\subseteq E_n$ for all  $n\in\omega$. The following conditions are equivalent:
\begin{enumerate}[(1)]
    \item $\mathcal{E}$ admits a stationary pseudointersection, that is, there is a stationary set $X$ such that, $X\setminus E_n$ is not stationary, for all $n\in\omega$. 
    \item $\displaystyle  \bigcap_{n\in\omega}E_n$ is stationary.
\end{enumerate}
\end{lemma}
\begin{proof}

$(1)\Rightarrow(2)$ Note that 
\[
X=(X\cap\bigcap_{n\in\omega}E_n)\cup(X\cap(\omega_1\setminus\bigcap_{n\in\omega}E_n))
\]
and one of the two sets forming the union must be stationary. On the other hand:
\[
X\cap(\omega_1\setminus\bigcap_{n\in\omega}E_n)=X\cap(\bigcup_{n\in\omega}\omega_1\setminus E_{n})=\bigcup_{n\in\omega}X\cap (\omega_1\setminus E_n)=\bigcup_{n\in\omega}X\setminus E_n,
\]
and as every $X\setminus E_n$ is not stationary, then neither is $\bigcup_{n\in\omega}X\setminus E_n$, i.e.  $ X\cap(\omega_1\setminus\bigcap_{n\in\omega}E_n)$ is not stationary. Necessarily $X\cap \bigcap_{n\in\omega}E_n$ is stationary and consequently $\bigcap_{n\in\omega}E_n$ also is stationary. 

$(2)\Rightarrow (1)$ In this case it is enough to take $X=\bigcap_{n\in\omega}E_n$. 
\end{proof}

\begin{corollary}
Let $\cI\subseteq \cP(\omega_1)$ be a $\cC$-independent family. The following conditions are equivalent:
\begin{enumerate}[(1)]
    \item $\cI$ is strongly $\cC$-independent.
    \item For every $f;\cI\rightarrow 2$ with $f$ countable, the collection $\{\cI^{f\restriction n}: n\in\omega\}$ admits a stationary pseudointersection\footnote{For $f\restriction n$ to make sense, it is enough to enumerate the domain of $f$ and so $f$ can be interpreted as a function in $2^\omega$.}.
\end{enumerate}
\end{corollary}

\begin{proposition}
A countable strongly $\cC$-independent family is not maximal, neither as a strongly $\cC$-independent family nor as a $\cC$-independent family.
\end{proposition}

\begin{proof}
Let $\cI=\{I_n: n\in\omega\}$ be a strongly $\cC$-independent family. For each $f\in2^\omega$ consider $X_f=\cI^f$. If $f\neq g$ then $X_f\cap X_g=\emptyset$. Now let $\{A_f,B_f\}$ be a partition of $X_f$ into stationary sets and define $A$ by:
\[
A=\bigcup_{f\in2^\omega} A_f.
\]
Let us see that $\cI\cup\{A\}$ is strongly $\cC$-independent. 
For this it is enough to see that for all $f\in2^\omega $, the sets $X_f \cap A$ and $X_f\setminus A$ are both stationary, however $ X_f\cap A=A_f$ and $X_f\setminus A=B_f$ 
are stationary sets.
\end{proof}

As we have seen, under $\mathsf{CH}$ there are countable $\cC$-independent families that are strongly $\cC$-independent, on the other hand (without extra hypothesis further than $\mathsf{ZFC}$) there are also countable $\cC$-independent families that are \textit{very far} from being strong. This means that there exists $\cI=\{I_n : n\in\omega\}\subseteq\cP(\omega_1)$ which is $\cC$-independent but that for every $h \in2^\omega$ such that $|h^{- 1}(\{0\})|=\omega $ we have that $\cI^h=\emptyset $; for example, to construct one of these families it is enough to take $\{X_n: n\in\omega\}$ a partition of $\omega_1$ into stationary sets and define $I_n$ as:
\[
I_n=\bigcup_{m\in C_n}X_m,
\]
where the $C_n$ are as in Example \ref{14}, in this way, the family $\{I_n : n\in\omega\}$ fulfills this property.

\section{Saturated ideals and $\cJ$-independent families}

Saturation of ideals has been closely related to the study of large cardinals, therefore it constitutes, as we will see in this section, a bridge between these cardinals and the existence of $\cJ$-independent families on them.

\begin{mydef}
Let $\cJ$ be an ideal on a cardinal $\kappa$. Then:
\begin{enumerate}
    \item $\cJ$ is $\lambda$-\emph{saturated} if for every collection $\{X_\alpha: \alpha\in\lambda\}\subseteq\cJ^+$ there exist $\beta<\gamma<\lambda$ such that $X_\beta\cap X_\gamma\in\cJ^+$.
    \item $\sat(\cJ)$ is the smallest $\lambda$ such that $\cJ$ is $\lambda$-saturated.
\end{enumerate}
\end{mydef}

\begin{lemma}\label{214}
Let $\cJ$ be an ideal on a cardinal $\kappa$ such that $\sat(\cJ)>\lambda$ for some cardinal $\lambda$. Then there exists a $\cJ$-independent family on $\kappa$ of cardinality $2^\lambda$.
\end{lemma}

\begin{proof}
Since $\cJ$ is not $\lambda$-saturated, it exists a collection $\left\{X_\beta:\beta\in\lambda\right\}\subseteq\cJ^+$ such that $X_\beta\cap X_\gamma\in\cJ$, if $\beta\neq\gamma$. Let $\cI=\{I_\alpha:\alpha\in 2^\lambda\}$ be an independent family of cardinality $2^\lambda$ on $\lambda$. For each $\alpha\in 2^\lambda$, let $\widehat{I_\alpha}\subseteq\kappa$ be defined as follows:
\[
\widehat{I_\alpha}=\bigcup\{X_\beta:\beta\in I_\alpha\}.
\]
Now set $\widehat{\cI}=\{\widehat{I_\alpha}:\alpha\in\kappa\}$. Clearly $\widehat{\cI}$ has cardinality $2^\lambda$, then the only thing left to prove is that it is an $\cJ$-independent family. 

Fix $s;2^\lambda\rightarrow2$, with $|s|<\omega$. We want to see that $\widehat{\cI}^s\in\cJ^+$. As $\cI$ is independent, it exists $\beta\in\cI^s$, but precisely the latter means that $X_\beta\subseteq\widehat{I_\alpha}$ for all $\alpha$ such that $s(\alpha)=0$ and $X_\beta\cap\widehat{I_\alpha}\in\cJ$ for all $\alpha$ such that $s(\alpha)=1$, i.e. $X_\beta\subseteq^*\widehat{\cI}^s$ and, since $X_\beta\in\cJ^+$, it follows that $\widehat{\cI}^s\in\cJ^+$.
\end{proof}

Next we will point out some relationships between the non-existence of strongly $\cJ$-independent families and the existence of large cardinals.

\begin{mydef}
If $\cJ$ is an ideal on $\kappa$, we say that $\cJ$ is $\kappa$-\emph{complete} if $\bigcup\mathcal{H}\in\cJ$, for every subfamily $\mathcal{H}\subseteq \cJ$ such that $|\mathcal{H}|<\kappa$.
\end{mydef}

\begin{theorem} (\cite{kanamori2008higher})\label{216}
Suposse that $\cJ$ is a $\kappa$-complete ideal on $\kappa$.
\begin{enumerate}
    \item (Tarski \cite{tarski1945ideale}) If $\cJ$ is $\lambda$-saturated with $2^{<\lambda}<\kappa$, then $\kappa$ is measurable.
    \item (Levy-Silver \cite{kanamori2008higher}) If $\cJ$ is $\kappa$-saturated and $\kappa$ is weakly compact, then $\kappa$ es measurable.
    \item (Kurepa \cite{kurepa1935ensembles}) If $\cJ$ is $\lambda$-saturated with $\lambda<\kappa$, then $\kappa$ has the tree property.
\end{enumerate}
\end{theorem}

\begin{corollary}
Suposse that $\cJ$ is a $\kappa$-complete ideal  on $\kappa$.
\begin{enumerate}
    \item If $\lambda<\kappa$, $2^{<\lambda}<\kappa$ and it does not exists a $\cJ$-independent family of cardinality $2^\lambda$, then $\kappa$ es measurable.
    \item If there is no $\cJ$-independent family of cardinality $2^\kappa$ and $\kappa$ es weakly compact, then $\kappa$ es measurable.
    \item If $\lambda<\kappa$ and there is no $\cJ$-independent family of cardinality $2^\lambda$, then $\kappa$ has the tree property.
\end{enumerate}
\end{corollary}

\begin{proof}
We will only prove the first part, the other two parts are analogous.

Since there is no $\cJ$-independent family of cardinality $2^\lambda$, then, by Lemma \ref{214}, we have that $\sat(\cJ)\leq \lambda$, i.e. $\cJ$ is $\lambda$-saturated, then by the first part of Theorem \ref{216} we have the desired result.
\end{proof}

Saturation of the ideal $\cJ$ is related to the existence of strongly $\cJ$-independent families.

\begin{proposition}\label{215}
Let $\cJ$ be an ideal on $\kappa$ and suppose that there exists a strongly $\cJ$-independent family of cardinality $\kappa$. Then $\sat(\cJ)\geq\kappa$. Furthermore, if $\kappa$ is regular then $\kappa$ is strongly inaccessible.
\end{proposition}

\begin{proof}
Let $\cI$ be a strongly $\cJ$-independent family of cardinality $\kappa$,  $\lambda<\kappa$ and $\cI_\lambda\subseteq\cI$ such that $|\cI_\lambda|=\lambda$. Then for every $h:\lambda\rightarrow 2$, we have that $\cI_\lambda^h\in\cJ^+$ and if $h\neq g$ then $\cI_\lambda^h\cap{\cI_\lambda}^g=\emptyset$, which proves that $\sat(\cJ)> 2^\lambda>\lambda$, and it finishes the proof. 
\end{proof}

The method in the previous proof has the advantage that it illustrates the fact that $\kappa$ is a strong limit cardinal, however the existence of a strongly $\cJ$-independent family of cardinality $ \kappa $ says even more about the saturation of $\cJ$: If $\cJ$ is an ideal on $\kappa$ and there is a strongly $\cJ$-independent family $\cI$ with cardinality $\kappa$, then $\sat(\cJ)> \kappa$.  Indeed, suppose that $\cI=\{X_\alpha:\alpha\in\kappa\}$ and for every $\beta\in\kappa$ let $Y_\beta=X_\beta\setminus \bigcup_{\alpha\in\beta}X_\alpha$. Note that, since $\cI$ is strongly $\cJ$-independent, $Y_\beta\in\cJ^+$, and if $\beta<\gamma<\kappa$ then $Y_\beta\cap Y_\gamma=\emptyset$. This proves that $\cJ$ is not $\kappa$-saturated (since $\{Y_\beta:\beta\in\kappa\}$ is a witness of that). 

\bibliography{ref}{}

\bibliographystyle{plain}

\end{document}